\documentclass{amsart}
\newtheorem{theorem}{Theorem}
\newtheorem{lemma}[theorem]{Lemma}
\theoremstyle{definition}
\newtheorem{definition}[theorem]{Definition}
\theoremstyle{remark}

\numberwithin{equation}{section}
\theoremstyle{corollary}

\theoremstyle{proposition}
\newtheorem{proposition}[theorem]{Proposition}

\newfont{\EUL}{eufm10 scaled 1000}

\newcommand\R{\mathbb{R}}

\newcommand\C{\mathbb{C}}

\newcommand\Z{\mathbb{Z}}

\newcommand\Ker{{\rm Ker}}

\newcommand\Lt{\mbox{\EUL t}}
\newcommand\h{\mbox{\EUL h}}
\renewcommand\k{\mbox{\EUL k}}

\newcommand\z{\mbox{\EUL z}}

\newcommand\Ad{{\rm Ad}}

\begin{document}
%
%
\title{Two-orbit K\"ahler manifolds and Morse Theory}
\author{Anna Gori and Fabio Podest\`a}
\address{Dipartimento di Matematica \lq U.Dini\rq\\ Viale Morgagni 67/A\\
50100 Firenze\\Italy} \email{gori@math.unifi.it}
\address{Dipartimento di Matematica e Appl.\ per l'Architettura\\
Piazza Ghiberti 27\\50100 Florence\\Italy}
\email{podesta@math.unifi.it} \subjclass{}
\begin{abstract} We deal with compact K\"ahler manifolds $M$ acted on by a compact Lie group $K$ of isometries, whose 
complexification $K^\C$ has exactly one open and one closed orbit in $M$. If the $K$-action is Hamiltonian, we obtain results on the cohomology and the $K$-equivariant cohomology of $M$.\end{abstract}
\maketitle
%
%
%
%
\section{Introduction}

Complete smooth complex algebraic varieties $M$ with an 
almost homogeneous action of a linear algebraic group $G$
have been studies by several authors (see e.g. \cite{Ak},
~\cite{Oe},~\cite{HS}). When the complement $N$ of the open $G$-orbit is homogeneous, the variety is called a two-orbit variety. 
Akhiezer (\cite{Ak}) gave a classification of two-orbit varieties 
in the case where $N$ is a complex hypersurface, while 
Feldm\"uller (\cite{F}) extended this classification under the hypeothesis that $N$ has codimension two and $G$ is reductive. 
Recently, Cupit-Foutou (\cite{C}) gave the full classification when $G$ is reductive.\\
We will be dealing with a compact K\"ahler manifold $M$, which
is acted on by a compact Lie group $K$ of isometries, whose 
complexification $K^\C$ acts on $M$ with two orbits. Under the 
hypothesis that the $K$-action is Hamiltonian with moment map $\mu$, we study the critical set of the squared moment map $||\mu||^2$,
applying standard Morse-theoretic results due to Kirwan 
(\cite{Kir1}) and obtaining information on the cohomology and 
$K$-equivariant cohomology of $M$.\\
Our main result is the following.
\begin{theorem}\label{t1}
Let $M$ be a compact K\"ahler manifold and let $K$ be a
compact connected Lie group of isometries acting on $M$ in a
Hamiltonian fashion. If the complexification $G$ of $K$ acts on
$M$ with two orbits, then
\begin{enumerate}
\item $K$ is semisimple and $M$ is simply connected projective
algebraic;
\item the Hodge numbers $h^{p,q}(M) = 0$ if $p\neq q$;
\item the function $f:M\to \R$ given by
$||\mu||^2$, where $\mu:M\to \k$ is the moment mapping, is a Bott-Morse function;
it has only two critical
submanifolds given by the closed $G$-orbit $N$, which realizes the
maximum of $f$ and by a $K$-orbit $S$, which realizes the minimum.
The Poincar\'{e} polynomial $P_M(t)$ of $M$ satisfies
$$P_M(t) = t^k\cdot P_N(t) + P_S(t) - (1 + t) R(t),$$
where $k$ is the real codimension of $N$ in $M$ and $R(t)$ is a
polynomial with positive integer coefficients. In particular 
$\chi(M) = \chi(N) + \chi(S)$;
\item if $\chi(M)> \chi(N)$, then $f$ is a perfect Bott-Morse function, i.e. 
$$P_M(t) = t^k\cdot P_N(t) + P_S(t).$$
\item The $K$-equivariant Poincar\'e series of $M$ is given by 
$$P_M^K(t) = t^k\cdot P_N^K(t) + P_S^K(t).$$

\end{enumerate}

\end{theorem}
We remark here that the statement of the Theorem holds for {\it
every\/} $K$-invariant K\"ahler form on $M$, not only for those
which are induced by the Fubini-Study K\"ahler form of a
projective space where the manifold might be $K^\C$-equivariantly
embedded.\\
In order to prove Theorem 1, we will need some properties of the
critical point set of the squared moment map. In \cite{Ne}, it is
proved that, given a compact Lie group acting holomorphically on
some $\Bbb P(V)$ with moment map $\mu$, then two critical points
of $||\mu||^2$ which belong to the same $K^\C$-orbit lie on the
same $K$-orbit; the proof of this fact is essentially based on the
explicit expression of the moment map for the
$K$-action on the projective space endowed with its canonical
Fubini-Study K\"ahler form. We slightly generalize such result to
a class of compact K\"ahler manifolds which are acted on by a
compact Lie group of isometries with moment map
$\mu$ and we put the following
\begin{definition} Given a compact K\"ahler manifold $(M,\omega)$ acted on by a
compact Lie group $K$ of isometries with moment map $\mu$, we say
that $\mu$ has the {\it Ness property\/} if two critical points
for $||\mu||^2$ which are in the same $G$-orbit belong to the same
$K$-orbit. \end{definition} Note that if $0\in \mu(M)$, then two
points in $\mu^{-1}(0)$ which lie in the same $G$-orbit belong to
the same $K$-orbit by a well-known fact due to \cite{Kir1},~p.97.
\par The next proposition gives sufficient conditions in order to
find invariant K\"ahler forms with the Ness property.
\begin{proposition}\label{prop} Let $M$ be a compact complex
manifold and let $K$ be a
connected, compact semisimple Lie group of holomorphic transformations of $M$ such
that $K^\C$ has an open orbit in $M$. If $\omega$ is a $K$-invariant
K\"ahler form, then there exists a
$K$-invariant K\"ahler form $\omega'$ which is cohomologous to
$\omega$ and whose corresponding moment map has the Ness property.
\end{proposition}
We will also need the following result, which might have an
independent interest (see also \cite{St}).
\begin{proposition}\label{max} Let $M$ be a compact K\"ahler manifold which is
acted on by a compact connected Lie group $K$ of isometries in a
Hamiltonian fashion with moment map $\mu$. If a point $x\in M$
realizes the maximum of $||\mu||^2$, then the orbit $K\cdot x$ is
complex.\end{proposition}
We remark here that Cupit-Foutou \cite{C} has deduced from her classification 
of two-orbit varieties that all such manifolds are spherical, i.e. a Borel subgroup of $G$ has an open orbit. Using the fact that the moment map relative to $K$ separates orbits when $M$ is spherical (see \cite{B}, \cite{HW}), we see that the critical sets $C_\beta$ of $f = ||\mu||^2$ (see section 2 for 
basic definitions) consist of single $K$-orbits, making the study of the Ness property unnecessary. We preferred to avoid using the fact that two-orbit 
varieties are spherical, with the hope of finding a new classification of 
such manifolds with a symplectic approach.  
%
%
\section{Basics on moment mappings and Morse theory}
We here briefly recall some basic results due to F.Kirwan
\cite{Kir1} about the geometry of the moment map; namely we will
first deal with a compact K\"ahler manifold $M$ which is acted on
effectively 
by a compact connected group of isometries $K$; we recall that $K$ acts automatically by holomorphic transformations and there is a holomorphic action 
on $M$ of the complex Lie group $G := K^\C$. \par
A smooth  map
$\mu:M\rightarrow \mathfrak{k}^*$ is said to be a moment map if
$$d\mu_x(X)(\xi)=\omega_x(X,\hat \xi)$$ for all $\xi\in \mathfrak
{k}$ and $X\in T_x M$,  where we denote by $\hat \xi$ the
fundamental Killing vector field on $M$ induced by the one
parameter subgroup $\exp(t\xi)$ for $\xi\in \mathfrak{k}$. If
$\mu$ is equivariant, then we will call $\mu$  an equivariant
moment map; it is known (see \cite{HW}) that an equivariant moment map exists 
if and only if $K$ acts trivially on the Albanese torus $Alb(M)$. 
In the sequel we will identify
$\mathfrak{k} \cong \mathfrak{k} ^*$ by means of a
$Ad(K)$-invariant scalar product $\langle,\rangle$ on the Lie
algebra $\mathfrak{k}$ and we will consider 
$\mu$ as a $\k$-valued map. \\ The critical points of the function
$f:M\to \R$ given by $f(x) := ||\mu||^2,$ where $\|\cdot\|$ is the
norm on
 $\mathfrak{k}$ induced by the scalar product  $\langle,\rangle$, are given
 by those points
$p\in M$ such that $\widehat{\mu(p)}_p = 0$; the
function $f$ is in general not a Bott-Morse function, while for
every $\beta\in \mathfrak{k}\setminus\{0\}$ the height function
$\mu_\beta:= \langle\mu,\beta\rangle$ is (see \cite{At}).
According to \cite{Kir1}, if $ \beta\in\mathfrak{k}$, $\beta\neq
0$, we denote by $Z_\beta$ the union of those connected components
of the critical set of $\mu_\beta$ on which $\mu_\beta$ assumes
the value $||\beta||^2$, while the Morse stratum $Y_\beta$
associated to $Z_\beta$ consists of all points of $M$ whose paths
of steepest descent under $\mu_\beta$ converge to points of
$Z_\beta$. The subset $Y_\beta$ is a complex submanifold (see
\cite{Kir1}, p. 89) and the Hessian $H(\mu_\beta)$ restricted to
$T_xY_\beta\times T_xY_\beta$ is positive semi-definite, where
$x\in Z_\beta$.\\ In \cite{Kir1} it is proved that $f$, although 
not nondegenerate in the sense of Bott, induces nevertheless a
smooth stratification $\{S_\beta|\;\beta\in{B}\}$ of $M$ for some
appropriate choice of a $K$-invariant metric on $M$. The stratum
to which a point of $M$ belongs is determined by the limit set of
its positive trajectory under the flow $-grad\;f$; the indexing
set $B$ is a finite subset of the positive Weyl chamber
$\mathfrak{t}_+$.\\ The strata $S_\beta$ are all locally closed
$K$-invariant submanifolds of $M$ of even dimension. In \cite{Kir1}
it is also proved that the set of critical points of $f$ is the
disjoint union of the closed subset $\{C_\beta=K(Z_\beta\cap
\mu^{-1}(\beta))\,|\;\beta\in{B}\}$ of $M$, and the image, under
$\mu$, of each connected component of the critical set of $f$ is a
single adjoint orbit in $\mathfrak{k}\simeq{\mathfrak{k}^*}$. More
precisely, for each $\beta\;\in{B}$, $C_\beta$ consists of those
critical points of $f$ whose image under $\mu$ lies in the adjoint
orbit of $\beta$; thus the function $f$ takes the constant value
$\|\beta\|^2$ on $C_\beta$.\\
 The introduction of the submanifolds $S_\beta$ is given in
\cite{Kir1} in the general setting of symplectic manifolds. If we
restrict our attention on K\"ahler ones, one can give  a different
characterization of the submanifolds $S_\beta$: \emph{a point $x$
belongs to $S_\beta$ if and only if $\beta$ is the closest point
to the origin of
$\mu(\overline{Gx})\cap{\mathfrak{t}_+}$} (\cite{Kir1},~p. 90).
%
%

\section{Proof of the main results}

%
%
We start proving Proposition 3, which will be useful in proving
Theorem 1. \par
\begin{proof}[Proof of Proposition~\ref{prop}] First of all we note that the
Albanese map $\alpha:M\to A(M):= Alb(M)$ is
$K$-equivariant and it induces a surjective homomorphism
$\alpha_*:G\to
Aut(A(M))\cong A(M)$ by \cite{Po}, where $G = K^\C$; since $G$ is
semisimple, we have that $A(M) = \{0\}$, i.e. the first Betti
number $b_1(M) = 0$. It then follows from \cite{Oe} that $M$,
being $G$-almost homogeneous, is simply connected algebraic and
$h^{2,0}(M) = 0$. This implies that $H^2(M,\C) =
H^1(M,\Omega_M^1)$ and therefore $H^2(M,\R) = H^{(1,1)}(M,\R)$,
where $H^{(1,1)}(M,\R) = H^2(M,\R)\cap H^1(M,\Omega_M^1)$. Since
$H^2(M,\Z)$ is a lattice of maximal rank in $H^2(M,\R)$, there
exists a basis $c_1,\ldots,c_m$ of $H^2(M,\R)$ with each $c_i\in
H^2(M,\Z)\cap K$, where $K$ denotes the open K\"ahler cone in
$H^2(M,\R)$. So, given a K\"ahler form $\omega$, we can find
positive real numbers $a_j$ such that $[\omega] =  \sum_j a_j
c_j$; since each $c_j$ is the Chern class of a positive line
bundle, this means that there exist embeddings $\phi_j:M\to {\Bbb
P}(V_j)$ such that $\omega$ is cohomologous to
$$\omega' := \sum_j \lambda_j\cdot (\phi_j)^*(\omega_{FS}({\Bbb P}(V_j))),$$
where $\omega_{FS}({\Bbb P}(V_j))$ denotes the K\"ahler form of
the Fubini-Study K\"ahler metric on ${\Bbb P}(V_j)$ and
$\lambda_j$ are real positive numbers. Note that $K$ acts
trivially on $H^1(M,{\mathcal O}^*)$ since $b_1(M) = 0$ implies
that $H^1(M,{\mathcal O}^*)\hookrightarrow H^2(M,\Z)$; this means
that the projective embeddings can be taken to be
$G$-equivariant.\par We will now prove that the moment map $\mu'$
corresponding to $\omega'$ has the Ness property, following the
lines of the proof in \cite{Ne},~p.1303; we will prove almost all
her statements  in full generality, while we will use the special
form of $\omega'$ only in Lemma \ref{l1}. Note that the
semisimplicity of $K$ implies the existence and the uniqueness of
$\mu'$. \par \vspace{0.3cm}
%
%
%
%
We now consider $x$ and $y$ two critical points that belong to the same
$G$-orbit. Recall that we can assume that $\mu'(x),\mu'(y)\neq 0.$
Firstly note that if $f'=\|\mu'\|^2$ then $f'(y)=f'(x).$ In fact
$x\in C_\beta\subset S_\beta$, where $\beta=\mu'(x)$ and, by the
$G$-invariance of $S_\beta$, $y=g\cdot x\in S_\beta$, since $y$ is
critical it must lie in $C_\beta$ on which $f'$ is constant. We
can also assume  that $\beta\in\mathfrak{t}_+$.  Let now take $P$
the parabolic subgroup corresponding to $\beta$; by the Levi
decomposition $P=L\cdot U$ where $L=K_\beta^\C$ is the Levi
subgroup of $P$, and $U$ its unipotent radical. We then have that
$G=KPU$, hence $y=g\cdot x=klu\cdot x$, for some $k\in K,l\in L$
and $u\in U$. Since $y$ is a critical point of $f'$ if and only if
the $K$-orbit $Ky$ consists of critical points of $f'$, it is
enough  to prove the result for $y=lu\cdot x$. Using this notation
we state
\begin{lemma}\label{l2} If  $\alpha\in \mathfrak{k}\setminus \mathfrak{k}_x$, 
then $\mu'_\alpha(\exp  ti\alpha\cdot x)$ is a strictly increasing
function of $t$.
\end{lemma}
\begin{lemma}\label{l1} Suppose $x$ and $y=lu\cdot x$ are critical points of $f'$ with $f'(x)=f'(y)\neq 0$. Then $u=e_G$,
 the identity in $G$.
\end{lemma}
Assuming the lemmas, we will deduce our claim. By Lemma \ref{l1},
$y=l\cdot x$; if $l\in G_x$ then $y=x$ and we are done. Otherwise
$y=l\cdot x$ and $l\notin G_x$. We will show that this implies
that $f'(y)>f'(x)$, a contradiction. Since $K_\beta^\C/K_\beta$ is
a symmetric space, there is a geodesic joining the cosets $[e]$
and $[l]$. Thus there is a one parameter subgroup $\exp
(ti\alpha_1)$, with $\alpha_1\in \mathfrak{k}_\beta$, such that
$l=\exp (t_0 i\alpha_1)k_1$ for $k_1\in K_\beta$ and $t_0\in
\R^+$. It suffices to prove the result for $y=\exp
(it\alpha_1)\cdot x$. Since $\beta\neq 0$ we take
$\alpha_0=\frac{\beta}{\|\beta\|}$. We can assume $\alpha_1$ to be
of length $1$ and perpendicular to $\alpha_0$, otherwise we
substitute $\alpha_1$ with its component in
$\mathfrak{k}_\beta\cap (\R\cdot \beta)^\perp$. Complete
$\alpha_0,\alpha_1$ to an orthonormal basis
$\{\alpha_0,\alpha_1\ldots\alpha_m\}$ of $\mathfrak{k}$. We have
$$f'(z)=\|\mu'(z)\|^2=\sum_{i=0}^m{ \mu'} _{\alpha_i}^2(z)$$ and
$$\|\mu'(y)\|^2=\|\mu'(\exp (ti\alpha_1)\cdot x)\|^2\geq {\mu'}_{\alpha_0}^2
(\exp ti\alpha_1\cdot x)+{\mu'}^2_{\alpha_1}(\exp ti\alpha_1\cdot
x).$$ One can easily prove that ${\mu'}_{\alpha_0}^2(\exp
ti\alpha_1\cdot x)$ is constant and equal to $\|{\mu'}(x)\|^2$; on
the other hand, $\mu'_{\alpha_1}(x)=0$ and ${\mu'}_{\alpha_1}(\exp
ti\alpha_1\cdot x)$ is strictly increasing,  by Lemma \ref{l2}, so
that ${\mu'}_{\alpha_1}^2(\exp ti\alpha_1\cdot x)$ is positive.
Hence $f'(y)>f'(x)$. This contradiction proves Proposition
\ref{prop}.
\end{proof}
Now we prove the two lemmas.
\begin{proof}[Proof of Lemma \ref{l2}] We will show that
 $\frac{d}{dt}\mu'_\alpha(\exp (it\alpha)\cdot x)>0$. Using the definition of $\mu'_\alpha$ we get
$$\frac{d}{dt}\mu'_\alpha(\exp (it\alpha)\cdot x)=\frac{d}{dt}<\mu'(\exp (it\alpha)\cdot x),\alpha>=
<d\mu'_{\exp ti\alpha\cdot x}(J\alpha),\alpha>$$ then, using the
definition of $\mu'$, we have that the last term is equal to $
\omega'_{\exp ti\alpha\cdot x}(J\hat\alpha,\hat\alpha)$ and
then equal to $g_{\exp
ti\alpha\cdot x}(\hat\alpha,\hat\alpha)$ which is strictly
positive.
\end{proof}
\begin{proof}[Proof of Lemma \ref{l1}]
Let $y=lu\cdot x$. Using the uniqueness of $\mu'$ we have that
$$\mu'(x)=\sum_i\lambda_i\mu'_i$$ where $\lambda_i$ are real
positive numbers, that have been already introduced, and $\mu'_i$
are the moment maps corresponding to the forms
$(\phi_i)^*(\omega_{FS}({\Bbb P}(V_i))).$ Using the same arguments as in 
\cite{Ne}, p.1305,
we have that, if $u\neq e$, for each $\mu'_i$ the following strict
inequality holds
\begin{equation}\label{ness}
{\mu'_i}_{\alpha_0}(lu\cdot
x)>{\mu'_i}_{\alpha_0}(x).
\end{equation}
We compute $\mu'_{\alpha_0}(lu\cdot x)$ and we prove that it is
positive: indeed
$$\mu'_{\alpha_0}(lu\cdot x)=
\sum_i \lambda_i<\mu'_i(lu\cdot x),\alpha_0>=\sum_i \lambda_i
{\mu'_i}_{\alpha_0}(lu\cdot x)$$ and, applying the strict
inequality (\ref{ness}), we have  that the last term is strictly
greater than $\sum_i \lambda_i
{\mu'_i}_{\alpha_0}(x)=\mu'_{\alpha_0}(x)=\|\mu'(x)\|>0$. Hence
$$f'(y)=\|\mu'(lu\cdot x)\|^2\geq {\mu'}^2_{\alpha_0}(lu\cdot x)>{\mu'}^2_{\alpha_0}(x)=\|\mu'(x)\|^2=f'(x)$$
so $u$ must be $e$ and Lemma \ref{l1} is proved.
\end{proof}
We next prove Proposition 4.
\begin{proof}[Proof of Proposition 4] Let $\beta = \mu(x)$, which
we can suppose to lie in $\Lt_+$. Then $x\in C_\beta =
K\cdot(Z_\beta\cap \mu^{-1}(\beta))$. We now claim that $Z_\beta =
\mu^{-1}(\beta)$. Indeed, if $p\in Z_\beta$, then $\mu_\beta(p) =
||\beta||^2$ and
$$||\beta||^2 \leq \langle\mu(p),\beta\rangle \leq ||\mu(p)||\cdot
||\beta|| \leq ||\beta||^2,$$
and therefore $\mu(p) = \beta$, i.e. $p\in \mu^{-1}(\beta)$.
Viceversa, if $p\in \mu^{-1}(\beta)$, then $||\mu(p)||^2$ is the
maximum value of $f:= ||\mu||^2$ and therefore $\hat{\beta}_p = 0$;
moreover $\mu_\beta(p) = ||\beta||^2$ and therefore $p\in
Z_\beta$. This implies that $\mu^{-1}(\beta)$ is a complex
submanifold and that $C_\beta = K\cdot \mu^{-1}(\beta)$. \\
We now claim that $S_\beta = C_\beta$. Indeed, if $\gamma_t(q)$ denotes the flow of
$-grad(f)$ through a point $q$ belonging to the stratum $S_\beta$,
then $\gamma_t(q)$ has a limit point in the critical subset
$C_\beta$; since $f(\gamma_t(q))$ is non-increasing for $t\geq 0$ and
$f(C_\beta)$ is the
maximum value of $f$, we see that $f(\gamma_t(q)) = ||\beta||^2$ for all
$t\geq 0$, that is $S_\beta\subseteq C_\beta$ and therefore
$S_\beta = C_\beta$.\\
This implies that $C_\beta = S_\beta$ is a smooth complex
submanifold of $M$ (see \cite{Kir1}, ~p. 89) and for every $y\in
\mu^{-1}(\beta)$, we have
$$T_yS_\beta = T_y(K\cdot y) + T_y(\mu^{-1}(\beta)).$$
Now, if $v\in T_y(\mu^{-1}(\beta))$, then $v = Jw$ for some $w\in
T_y(\mu^{-1}(\beta))$ and for every $X\in \k$ we have
$$0 = \langle d\mu_y(w),X\rangle = \omega_y(w,{\hat X}_y) =
\omega_y(Jv,{\hat X}_y) = g_y(v,{\hat X}_y),$$ meaning that
$T_y(\mu^{-1}(\beta))$ is $g$-orthogonal to $T_y(K\cdot y)$. Since
both $S_\beta$ and $\mu^{-1}(\beta)$ are complex, this implies
that $K\cdot y$ is a complex orbit.\end{proof}

%
%

\begin{proof}[Proof of Theorem 1]
We will first prove that $K$ is semisimple. Let $\z$ be the Lie
algebra of the center of $\k$, the Lie algebra of $K$, and let
$X\in \z$ be a vector such that $\exp(X)$ generates the torus
$Z^o(K)$. We will show that $X$ acts trivially on $M$; let $C$ denote the zero set of $X$.\par
If $\mu:M\to \k$ denotes a moment mapping for the
$K$-action, then $\phi:\pi\circ\mu:M\to \R\cdot X$, where
$\pi:\k\to \R\cdot X$ denotes the orthogonal projection w.r.t. an
$\Ad(K)$-invariant scalar product on $\k$, is a function such that
$d\phi_x(Y) = \omega(Y,X)$ for all $Y\in T_xM$ and all $x\in M$,
where $\omega$ is the K\"ahler form of $M$. Since critical points
of $\phi$ are exactly the zeros of $X$ in $M$, we have that there
are at least two connected components of the zero set of $X$; on
the other hand if $x\in \Omega$ is a point where $X$ vanishes and
$\Omega = G/H$ is the open $G$-orbit, then $\z^\C\subset \h$ and
therefore $Z^o(K)$ acts trivially on $\Omega$, hence on the whole
$M$. But then the zero set of $X$ lies in the closed orbit $N$,
which is connected, so that $C = N$, a contradiction.\\
Arguing as in the proof of Proposition 3, we see that the first
Betti number $b_1(M) = 0$. It then follows from \cite{Oe} that
$M$, being $G$-almost homogeneous, is simply connected and
projective algebraic.\par We now consider the function $f:=
||\mu||^2$ and its critical set $C$. We first note that $C$ does
not coincide with $M$, i.e. the function $f$ is not constant;
indeed otherwise, by Proposition 4, every $K$-orbit would be
complex, contrary to the uniqueness of the closed $G$-orbit.\par

By Proposition 4, we see that the closed orbit $N$ consists
exactly of those points at which the function $f$ realizes the
maximum value; we denote by $\beta_{\text{max}}$ the point
$\mu(N)\cap\Lt_+$. \par
 We now denote by $\beta_{\text{min}}$ the point in
the convex polytope $\mu(M)\cap \Lt_+$ which is the closest point
to the origin; then $||\beta_{\text{min}}||^2$ is the minimum
value of $f$. Moreover we know (\cite{Kir1}) that a point $x$
belongs to the stratum $S_\beta$ for some $\beta\in \Lt_+$ if and
only if $\beta$ is the closest point to the origin of
$\mu(\overline{G\cdot x})\cap \Lt_+$; it then follows that the
whole orbit $G\cdot x$ belongs to $S_{\beta_{\text{min}}}$ and
therefore there are exactly two critical values of $f$, that is
$||\beta_{\text{max}}||^2$ and $||\beta_{\text{min}}||^2$; moreover
$S_{\beta_{\text{min}}} = \Omega$, where $\Omega$ is the open $G$-orbit. \par
Since $M = S_{\beta_{\text{max}}} \cup S_{\beta_{\text{min}}}$, we
see that a critical point of $f$ must lie either in the closed
orbit $N$ or in $C_{\beta_{\text{min}}}$. We here prove that
$C_{\beta_{\text{min}}}$ consists of a single $K$-orbit. Indeed we have
seen in Proposition \ref{prop} that there exists a K\"ahler form
$\omega'$ on $M$ which is cohomologous to $\omega$ and whose
corresponding moment map has the Ness property. 
Since $\omega_t = t\omega' + (1 - t)\omega$, $t\in [0,1]$, is a 
homotopy by symplectic forms connecting $\omega$ with $\omega'$, by the
equivariant Moser homotopy Theorem (\cite{SM},~p. 91) there exists a 
$K$-equivariant symplectomorphism 
$$F:(M,\omega)\to (M,\omega').$$ 
Now note that the (unique) moment map corresponding to $\omega$ is 
$\mu=\mu'\circ F$, hence the
critical points of $f$ are taken by $F$ to critical points of
$f'$. Using the same arguments as above, we see that 
$\Omega = S'_{\beta_{\text{min}}}$, where $S'$ denotes the stratum w.r.t. 
$\mu'$;  hence the Ness property implies that
the minimal critical set for $f'$ consists of a single $K$-orbit and,
using the $K$-equivariance of $F$, we conclude that also
$C_{\beta_{\text{min}}}$ consists of a single $K$-orbit. We  have
thus proved the following
\begin{proposition} The inverse images $\mu^{-1}(K\cdot
\beta_{\text{max}})$ and $\mu^{-1}(K\cdot\beta_{\text{min}})$
consist exactly of two $K$-orbits and the first one is the closed
$G$-orbit.\end{proposition}
We will call $S$ the $K$-orbit such that $\mu(S) =
K\cdot\beta_{\text{min}}$. We now consider a maximal torus $T$ inside
$K$ and prove the following easy Lemma.
\begin{lemma} The fixed point set $M^T$ is contained in the
critical set of $f$.\end{lemma}
\begin{proof} Indeed if $x\in M^T$, then $T\subseteq K_x$ and by
the $K$-equivariance of $\mu$, we have that $\Ad(T)\mu(x) =
\mu(x)$. Since $T$ is maximal, we have that $\mu(x)\in \Lt$, where
$\Lt$ is the Lie algebra of $T$, and therefore $\mu(x)\in
\k_x$\end{proof}
From this we see that the set $M^T$ is finite, because a maximal
torus of $K$ fixes at most a finite number of points in a
$K$-homogeneous space. Using \cite{CL} we deduce that the Hodge
numbers $h^{p,q} = 0$ for $p\neq q$.\\
We have now to prove that the function $f$ is a Bott-Morse
function and the rest of the statement (iii) will follow from
standard Morse theory.\\
Let $Z_\beta^{min}$ be the subset of $Z_\beta$ consisting of those points $x\in Z_\beta$ such that the limit points of the path of steepest descent from $x$ for the function $\|\mu-\beta\|^2$ on $Z_\beta$ lie in $Z_\beta\cap \mu^{-1}(\beta)$, and let $Y_\beta^{min}$ be the inverse inverse image of $Z_\beta^{\min}$
 under the retraction $p_\beta:Y_\beta\rightarrow Z_\beta$. Then $Y_\beta^{min}$ is an open subset of $Y_\beta$ and retracts to $Z_\beta^{min}.$\\
With this notation we recall (see \cite{Kir1}) that the
$G$-open orbit $\Omega$ has the structure of a fiber bundle
$G\times_{P_\beta}Y_\beta^{min}$, where  $P_\beta$
is the parabolic subgroup of $G$ such that $G/P_\beta =
K/K_\beta$ and we recall also that $T_pY_\beta \cap T_x (K\cdot p) =
T_p(K_\beta\cdot p)$ ( see Lemma 4.10, p. 48 in \cite{Kir1}).\\
We start showing that the Hessian of $f$ at a
critical point $p$ belonging to the minimum orbit $S$ is
non degenerate; we can suppose that $\mu(p) = \beta$, where we have
put $\beta = \beta_{\text{min}}$.
We have
$$\dim K/K_p + \dim Y_\beta - \dim(T_pY_\beta\cap T_p(K\cdot p))=$$
$$=\dim K/K_p + \dim Y_\beta - \dim (K_\beta/K_p) = \dim K/K_\beta
+ \dim Y_\beta = \dim M,$$
hence
$$T_p(K\cdot p) + T_p{Y_\beta} = T_pM.$$
If we consider $W_p :=
{(T_p (K_\beta\cdot p))}^\perp \cap T_p Y_\beta$, then $W_p$ is a
complement of $T_p(K\cdot p)$ in $T_pM$;
we will compute the Hessian of $f$ on vectors $X\in {W_p}$, namely we will show that
$H(f)_p(X,X)> 0$ if $X\in {W_p}$, $X\neq 0$. \\
We fix an orthonormal basis $\{X_1,\dots,X_l\}$ of the Lie algebra $\k$ w.r.t. the invariant
scalar product $\langle,\rangle$; then we can write for all $x\in M$,
$\mu(x) = \sum_i\mu_i(x)X_i$. Moreover we denote by ${\mathcal M}$
the vector field on $M$ given by $\mathcal{M}_x = \widehat{\mu(x)}$ for $x\in M$, where
for any $X\in \k$ we denote by $\hat X$ the induced Killing vector field on $M$.
Then we have
$$df_x(X) = 2\langle d\mu(x)(X),\mu(x)\rangle = 2\omega_x(X,{\mathcal M}_x).$$
It follows that $$H(f)_x(X,X) = \omega_x(X,\nabla_X{\mathcal M}),$$
where $\nabla$ denotes the Levi Civita connection of the K\"ahler metric $g$.
We then have that
$$H(f)_p(X,X) = 2\sum_i\omega_p(X,X(\mu_i)\hat{X_i} +
\mu_i(p)\nabla_X\hat{X_i}).$$
Now we observe that
$$\omega_p(X,\hat{X_i}) = \langle d\mu_x(X),X_i\rangle = X(\mu_i)$$
and therefore
\begin{displaymath}
H(f)_p(X,X) = 2\sum_i\omega_p(X,\hat{X_i})^2 + \omega_p(X,\nabla_X\hat{\beta}).
\end{displaymath}
If $\mu_\beta$ denotes as usual the height function $\langle\mu,\beta\rangle$, we see
that $H(\mu_\beta)_p(X,X) = 2\omega_p(X,\nabla_X\hat\beta)$,
so that
\begin{displaymath}
H(f)_p(X,X) = 2\sum_i\omega_p(X,\hat{X_i})^2 + 2H(\mu_\beta)_p(X,X).
\end{displaymath}
We now distinguish two cases, whether $\beta = 0$ or $\beta\neq 0$.\\
If $\beta = 0$, then $H(F)_p(X,X) = 2\sum_i\omega_p(X,\hat{X_i})^2$ for all $X\in T_p(K\cdot p)^\perp$
and therefore $X$ belongs to the nullity of $H(f)_p$ if and only
if $JX$ is orthogonal to the tangent space $T_p(K\cdot p)$. But
$\beta = 0$ implies that $J(T_p(K\cdot p))$ is contained in the
normal space $T_p(K\cdot p)^\perp$; this together with the fact
that $p$ belongs to the open $K^\C$-orbit implies that
$J(T_p(K\cdot p)) = T_p(K\cdot p)^\perp$ and therefore $X = 0$.
This shows that $H(f)_p$ is non degenerate.\\
If $\beta\neq 0$, then we choose $X\in W_p.$
Denote by $V$ the tangent space $T_p(K\cdot p)$. Note that
$\omega_p(X,\hat{X}_i)=0$ if and only if $d\mu_p(X)=0$ i.e. if and only if $X\in
(JV)^\perp$ and that $H({\mu_{\beta}})_p (X,X)=0$ if and only if $X\in T_p Z_\beta$.
We shall prove that, if $H(f)_p(X,X)=0$, then $X$ must be orthogonal to $V$; hence,
since $V+JV=T_p M$, $X$ must be $0$. In fact, recall that $p$ is a critical point
for $\mu_\beta$ and we have $\hat{\beta}_p=0$. We consider the skew operator
$\nabla\hat{\beta}:T_p M\rightarrow T_p M$; since $e^{t{\nabla\hat{\beta}}} =
d(exp(t\beta))_p$, we have that
$\nabla\hat{\beta}$ maps $V$ into itself and $T_p Z_\beta=Ker \nabla\hat{\beta}$.
The tangent space $V$ then  splits as the direct sum of $\Ker (\nabla\hat{\beta}_{|_V})$ and
its orthogonal complement $V'$.
Note that $\nabla\hat{\beta}_{|_{V'}}$ is surjective on $V'$, hence, if $Y\in V'$
we have $Y=\nabla\hat{\beta}(Y')$ for some $Y'\in V'$.
Using that $\nabla\hat{\beta}$ is skew, we argue that $X\perp V'$. Finally note that
$\Ker (\nabla\hat{\beta}_{|_V}) = T_p (K_{\beta}\cdot p)$, hence, since $X\in W_p$,
it is orthogonal to $T_p (K_{\beta}\cdot p)$ and we get our claim.
\end{proof}
The statements in (iii) and (v) follow from the general theory in \cite{Kir1}.
As for (iv), we note that if $\chi(M) > \chi(N)$, then $S$ is a compact homogeneous space of positive Euler characteristic; therefore all odd Betti numbers 
of $S$ vanish (see \cite{BH}) and therefore $f$ is perfect.

%
%

\end{document}